\documentclass[12pt]{article}
\usepackage[latin1]{inputenc}
\usepackage{amssymb,amsmath,amsfonts}
\newtheorem{theorem}{Theorem}
\newtheorem{proposition}[theorem]{Proposition}
\newtheorem{lemma}[theorem]{Lemma}

\newtheorem{corollary}[theorem]{Corollary}

\newcommand{\R}{\mathbb{R}}

\newcommand{\Q}{\mathbb{Q}}
\newcommand{\Sf}{\mathbb{S}}
\newcommand{\C}{\mathbb{C}}

\newcommand{\spa}{\mbox{span}}

\newcommand{\po}{{\hspace*{-1ex}}{\bf .  }}
\newcommand{\ii}{isometric immersion }

\newcommand{\nab}{\tilde\nabla}

\def\D{{\cal D}}
\def\<{\langle}
\def\n{\nabla}

\def\>{\rangle}
\def\a{\alpha}

\def\bea{\begin{eqnarray*} }
\def\eea{\end{eqnarray*} }
\def\be{\begin{equation} }
\def\ee{\end{equation} }

\def\proof{\noindent{\it Proof:  }}
\def\qed{\ifhmode\unskip\nobreak\fi\ifmmode\ifinner
\else\hskip5 pt \fi\fi\hbox{\hskip5 pt \vrule width4 pt
height6 pt  depth1.5 pt \hskip 1pt }}
\begin{document}

\title{Complete minimal submanifolds with nullity\\ in Euclidean 
spheres}
\author{M.\ Dajczer, Th.\ Kasioumis,
A.\ Savas-Halilaj and Th.\ Vlachos}
\date{}
\maketitle

\renewcommand{\thefootnote}{\fnsymbol{footnote}} 
\footnotetext{The third author would like to acknowledge 
financial support from the grant DFG SM 78/6-1.}  
\renewcommand{\thefootnote}{\arabic{footnote}}

\begin{abstract} 
In this paper we investigate  $m$-dimensional complete minimal submanifolds 
in Euclidean spheres with  index of relative nullity  at least $m-2$ at any point.
These are austere submanifolds in the sense of Harvey and Lawson \cite{harvey}
and were initially studied by Bryant \cite{br}.
For any dimension and codimension there is an abundance of non-complete examples 
fully described by Dajczer and Florit  \cite{DF2} in terms of a  
class  of surfaces, called elliptic, for which the ellipse of curvature of
a certain order is a circle at any point. 
Under the assumption of completeness, it turns out 
that any submanifold is either totally geodesic or has dimension three. 
In the latter case there are plenty of examples, even compact ones. 
Under the mild assumption that the Omori-Yau maximum principle  
holds on the manifold, a trivial condition in the compact case, 
we provide a complete local parametric description of the submanifolds 
in terms of $1$-isotropic surfaces in  Euclidean space. These are the minimal 
surfaces for which the standard ellipse of curvature is a circle
at any point. For these surfaces, there exists a Weierstrass type representation 
that generates all simply connected ones.
\end{abstract} 

Let $M^m$ be a complete $m$-dimensional Riemannian manifold. In \cite{dksv} 
we considered the case of minimal isometric immersions  into Euclidean space 
$f\colon M^m\to\R^n$, $m\geq 3$,
satisfying that the index of relative nullity  is at least $m-2$ 
at any point. Under the mild assumption that the Omori-Yau maximum principle 
holds on $M^m$, we concluded that any $f$ must be ``trivial", namely, just a 
cylinder over a complete minimal surface. This result is global in nature 
since for any dimension there are plenty of non-complete examples  other 
than open subsets of cylinders.

It is natural to expect rather different type of conclusions when considering 
a similar global problem for minimal isometric immersions into nonflat space 
forms. For instance,  for submanifolds in the hyperbolic space one would 
guess that under the same condition on the relative nullity index there 
exist many non-trivial examples, and that a kind of triviality  conclusion 
will only hold under a strong additional assumption.

This paper is devoted to the case of minimal submanifolds in round spheres. 
Thus, in the sequel  $f\colon M^m\to\Sf^n$, $m\geq 3$, will be a minimal isometric
immersion into Euclidean sphere with index of relative nullity at  least $m-2$ 
at any point. As in the Euclidean case, it is known that are plenty  of 
non-complete examples of any dimension.      
On one hand, under the assumption of  completeness of $M^m$ we have that $f$ 
has to  be totally geodesic unless $m=3$. On the other hand, for dimension $m=3$ 
we will see that there are plenty of non trivial examples, even compact ones.
\medskip

Notice that the class of submanifolds studied in this paper are precisely the 
minimal $\delta(2)$-ideal submanifolds considered by Chen; see \cite{byc}.
\medskip

From the results in \cite{DF2} it follows that any example for dimension $m=3$ 
can locally be constructed  as follows:
\medskip

Let $g\colon L^2\to\R^{n+1}$, $n\ge 4$, be an elliptic surface whose first 
curvature ellipse is always a circle.  Then, the map $\psi_g\colon T^1L\to\Sf^n$ 
defined on the  unit tangent bundle of $L^2$  and given by  
\be\label{psi}
\psi_g(x,w)=g_*w
\ee
parametrizes (outside singular points) a minimal immersion $f\colon M^3\to\Sf^n$  
with index of relative nullity at least one at any point.
\vspace{1ex}

Minimal surfaces are elliptic, but the latter class of surfaces is much larger. 
In fact, that a surface $g\colon L^2\to\R^n$ is elliptic means that 
given a (hence any) basis $X,Y$ of the tangent plane $T_xL$ at any $x\in L^2$ 
the second fundamental form $\alpha_g\colon TL\times TL\to N_gL$ of $g$ 
satisfies
$$
a\alpha_g(X,X)+2b\alpha_g(X,Y)+c\alpha_g(Y,Y)=0
$$
where $a,b,c\in\R$ verify $ac-b^2>0$. Equivalently, in any local system of 
coordinates $(u,v)$ of $L^2$ any coordinate function of $g$ is a solution of the 
(same) elliptic PDE of type
$$
a\frac{\partial^2}{\partial u^2}+2b\frac{\partial^2}{\partial u\partial v}
+c\frac{\partial^2}{\partial v^2}+d\frac{\partial}{\partial u}
+e\frac{\partial}{\partial v}=0
$$
where $a,b,c,d,e$ are smooth functions  such that  $ac-b^2>0$.

Our main result shows that for  complete submanifolds the condition 
of ellipticity of the generating surface has to be restricted to minimality.
For minimal surfaces, the first curvature ellipse as an elliptic surface 
coincides with the standard ellipse of curvature, namely, the image in the 
normal space of the second fundamental form restricted to the unit 
circle in the tangent plane.

\begin{theorem}\label{main}\po
Let  $f\colon M^m\to\Sf^n$, $m\geq 3$, be a  minimal isometric immersion with 
index of relative nullity at least $m-2$ at any point. If $M^m$ is complete 
then $f$ is totally geodesic unless $m=3$.  Moreover, if the Omori-Yau 
maximum principle holds on $M^3$ then along an open dense subset  $f$ is 
locally parametrized  by (\ref{psi}) where $g\colon L^2\to\R^{n+1}$ is 
a minimal surface whose first curvature ellipse is always a circle.
\end{theorem}

A minimal surface $g\colon L^2\to\R^n$ whose first curvature ellipse 
is a circle at any point is called a \mbox{$1$-isotropic} surface. 
The above result should be complemented by the fact that there is the
Weierstrass type representation from  \cite{dg2} that generates all 
simply connected \mbox{$1$-isotropic} surfaces.  It goes as follows:
start with any nonzero holomorphic map in a simply connected domain 
$\alpha_0\colon U\subset\C\to\C^{n-4}$ and let  $\alpha_1\colon U\to\C^{n-2}$
be given by 
$$
\alpha_1=\beta_1\left(1-\phi_0^2,i(1+\phi_0^2),2\phi_0\right)
$$
where $\phi_0=\int^z\alpha_0dz$ and $\beta_1\neq 0$ is any holomorphic function. 
Now let $\alpha_2\colon U\to\C^n$ be given by
$$
\alpha_2=\beta_2\left(1-\phi_1^2,i(1+\phi_1^2),2\phi_1\right)
$$
where $\phi_1=\int^z\alpha_1dz$ and $\beta_2\neq 0$ is any holomorphic function. 
Then $g=\mbox{Re}\{\alpha_2\}$ is a $1$-isotropic surface in $\R^n$. 

\vspace{1ex}

The Omori-Yau maximum principle  holds on $M^m$
if for any function $\varphi\in C^2(M)$ bounded from above 
there exists a sequence of points $\{x_j\}_{j\in\mathbb{N}}$
such that
$$
\varphi(x_j)>\sup_M\varphi-1/j,\quad \|\nabla\varphi(x_j)\|\le 1/j\quad
\text{and}\quad\Delta\varphi(x_j)\le 1/j
$$
for any $j\in\mathbb{N}$.
There are fairly general assumptions that imply the validity of the 
Omori-Yau maximum principle on a Riemannian manifold; see \cite{amr}. 
For instance, it is applicable on complete Riemannian manifolds whose 
Ricci curvature does not decay fast to $-\infty$.
\vspace{1ex}

We see next that there are plenty of complete (even compact) examples of three-dimensional 
minimal submanifolds in spheres with index of relative nullity at least one
at any point.
\vspace{1ex}

\noindent\emph{Hopf lifts:}  If $g\colon L^2\to\mathbb{CP}^n$, $n\ge 2$, is
a substantial holomorphic curve, then the Hopf fibration 
$\mathcal{H}\colon\Sf^{2n+1}\to\mathbb{CP}^n$ induces a circle bundle 
$M^3$ over $L^2$. This lifting induces an immersion
$f\colon M^3\to\Sf^{2n+1}$ such that $g\circ\pi=\mathcal{H}\circ f$,
where $\pi\colon M^3\to L^2$ is the projection map.
Such submanifolds are minimal with 
index of relative nullity at least one if $n=2$ (see \cite{dillen}) or 
$n=3$ (see \cite{lotay}). 
Moreover, if $L^2$  is compact, then $M^3$ is also  compact.
\vspace{1ex}

\noindent \emph{Tubes over minimal 2-spheres}:
Due to the work of Calabi, Chern, Barbosa and others,  it is known 
that minimal $2$-spheres in spheres are  pseudoholomorphic (isotropic) 
in substantial even codimension.  Calabi \cite{calabi} proved that any such 
surface in $\Sf^{2n}$ is nicely curved if its area is $2\pi n(n+1)$, and 
Barbosa showed \cite{barbosa} that the space of  these surfaces is diffeomorphic 
to $SO(2n+1,\C)/SO(2n+1,\R)$.  According to Proposition \ref{polar-cl} below 
such surfaces produce examples of compact three-dimensional minimal submanifolds 
in $\Sf^{2n}$ with index of relative nullity one. 
\vspace{1ex}

Among the second family of examples given above, there are the submanifolds 
produced from pseudoholomorphic surfaces $g\colon\Sf^2\to\Sf^6$ with area 
$24\pi$ which are holomorphic with respect to the nearly Kaehler structure 
in $\Sf^6$. For instance, this is the situation of the Veronese surface in $\Sf^6$. 
In this case, the compact submanifolds $M^3$ are  
Lagrangian (also called totally real) in $\Sf^6$; see \cite{dillen}.

\begin{corollary}\label{Lag}\po
Let  $f\colon M^3\to\Sf^6$ be an isometric immersion with index 
of relative nullity at least one at any point. Assume that $f$ is Lagrangian 
with respect to the nearly Kaehler structure in $\Sf^6$.
If $M^3$ is complete and the Omori-Yau maximum principle  holds, then $f$ 
is locally parametrized by (\ref{psi}) along an open dense subset of $M^3$ 
where $g$ is a $2$-isotropic surface in $\R^6$ (respectively, $\R^7$) and 
$f$ is substantial in $\Sf^5$ (respectively, $\Sf^6$). 
\end{corollary}

That the surface $g$ is $2$-isotropic means that it is $1$-isotropic 
and that the second ellipse of curvature is also a circle at any point.  
Hence, in the case of $\R^6$ we have that $g$ is congruent to a 
holomorphic curve in $\C^3\equiv\R^6$.  
\vspace{1ex}

It follows from the results in \cite{DV0} that 
the universal cover of any of the complete three-dimensional submanifolds  
considered in Theorem \ref{main} admits a one-parameter associated family 
of isometric immersions of the same type.  Moreover, that family is trivial 
if and only if the (local) generating minimal surface is congruent to a 
holomorphic curve. We refer to Lotay \cite{lotay} for a discussion about 
the existence of such an associated family in the case of yet another 
family of examples.

\section{The relative nullity foliation}

In this section, we recall basic facts on the relative nullity foliation 
for submanifolds in space forms to be used in the sequel without 
further reference.\vspace{1ex}

Let $f\colon M^m\to\Q_c^n$ denote an \ii of an $m$-dimensional Riemannian 
manifold $M^m$ into the Euclidean space $\R^n$ ($c=0$) or the unit sphere 
$\Sf^n$ ($c=1$). 
\vspace{1ex}

The \emph{relative nullity subspace} of $f$ at $x\in M^m$ is the tangent subspace
given by
$$
\D(x)=\{X\in T_xM\colon \alpha_f(X,Y)=0\;\;\text{for all}\;\;Y\in T_xM\}
$$
where $\alpha_f\colon TM\times TM\to N_fM$ denotes the second fundamental form
of $f$.
The dimension $\nu(x)$ of $\D(x)$ is called the \emph{index of relative nullity} 
of $f$ at $x\in M^m$. 
\vspace{1ex}

For simplicity, in the sequel we call $k(x)=m-\nu(x)$ the \emph{rank} of $f$ 
at $x\in M^m$. Notice that $k(x)$ is the rank of the Gauss map of $f$ at $x\in M^m$.
If $f$ has constant rank on an open subset $U\subset M^m$, it is an elementary fact 
that the relative nullity distribution $\D$ along $U$ is integrable and that the 
leaves of the \emph{relative nullity foliation} are  totally geodesic submanifolds of 
$M^m$  and under $f$ of $\Q_c^n$. 
\vspace{1ex}

Let  $U\subset M^m$ be an open subset where the index of relative nullity 
$\nu=s>0$  is  constant. The following is a well known fundamental result 
in the theory of isometric immersions (cf.\ \cite{da}).

\begin{proposition}\label{comp}\po Let $\gamma\colon [0,b]\to M^m$ be a 
geodesic curve such that $\gamma([0,b))$ is contained in a leaf of relative nullity 
contained in $U$. Then also $\nu(\gamma(b))=s$.
\end{proposition}

In the sequel, we decompose any tangent vector field $X$ on $M^m$ as $X=X^v+X^h$ 
according to the orthogonal splitting $TM=\D\oplus \D^{\perp}$.
\vspace{1ex}

The \textit{splitting tensor} $C\colon\D\times\D^{\perp}\to\D^{\perp}$ 
is given by
$$
C(T,X)=-{\nabla}^{h}_XT=-(\nabla_XT)^h
$$
for any $T\in\D$ and $X\in\D^{\perp}$. We also regard $C$ as a map
$C\colon\Gamma(\D)\to\Gamma(\operatorname{End}(\D^{\perp}))$.
\vspace{1ex}

The following differential equations for the endomorphism $C_T=C(T,\cdot)$
are a well known easy consequence of the Codazzi equation:
\be\label{C1}
\n_S C_T=C_T C_S+C_{\n_ST}+c\<T,S\>I
\ee
where $I$ is the identity map, and
\be\label{C2}
(\nabla^{h}_{X}C_T)Y-(\nabla^{h}_{Y}C_T)X
=C_{\nabla^{v}_{X}T}Y-C_{\nabla^{v}_{Y}T}X
\ee
for any $S,T\in\Gamma(\D)$ and $X,Y\in\Gamma(\D^{\perp})$. For a proof we 
refer to \cite{da} or \cite{dg}.

\section{Elliptic submanifolds}

In this section, we recall from \cite{DF2} the notion of elliptic submanifold 
of a space form as well as several of their basic properties.
\vspace{1ex}

Let $f\colon M^m\to\Q_c^n$ be an isometric immersion.
The $\ell^{th}$\emph{-normal space} $N^f_\ell(x)$  of $f$
at $x\in M^m$ for $\ell\ge 1$ is defined as
$$
N^f_\ell(x)=\spa\big\{\alpha_f^{\ell+1}(X_1,\ldots,X_{\ell+1}):
X_1,\ldots,X_{\ell+1}\in T_xM\big\}.
$$
Here $\alpha_f^2=\alpha_f$ and for $s\geq 3$ the so called 
$s^{th}$\emph{-fundamental form} is the symmetric tensor 
$\alpha_f^s\colon TM\times\cdots\times TM\to N_fM$  defined inductively by
$$
\alpha_f^s(X_1,\ldots,X_s)=\pi^{s-1}\left(\nabla^\perp_{X_s}\cdots
\nabla^\perp_{X_3}\alpha_f(X_2,X_1)\right)
$$
where $\pi^k$ stands for the projection onto 
$(N_1^f\oplus\cdots\oplus N_{k-1}^f )^{\perp}$.

An isometric immersion $f\colon M^m\to\Q_c^n$ of rank $2$ is called \emph{elliptic}  
if there exists a (necessary unique up to a sign) almost complex 
structure $J\colon\D^\perp\to \D^\perp$ such that  the second fundamental 
form satisfies
$$
\alpha_f(X,X)+\alpha_f(JX,JX)=0
$$
for all $X\in\D^\perp$. Notice that $J$ is orthogonal if and only $f$ is 
minimal. 

Let $f\colon M^m\to\Q_c^n$ be substantial and elliptic. 
The former means that the codimension cannot be reduced.
Assume also that $f$ is \emph{nicely curved} which means that for any $\ell\geq 1$
all subspaces $N^f_\ell(x)$ have constant dimension and thus form subbundles
of the normal bundle. Notice that any $f$ is nicely curved along connected 
components of an open dense subset of $M^m$.
Then, along that subset the normal bundle splits orthogonally and smoothly as
\be\label{splits}
N_fM=N^f_1\oplus \cdots \oplus N^f_{\tau_f}
\ee
where all $N^f_\ell$'s have rank $2$, except possibly the last one 
that has rank $1$ in case the codimension is odd.
Thus, the induced bundle $f^*T\Q_c^n$ splits as 
$$
f^*T\Q_c^n=f_*\D\oplus N^f_0\oplus N^f_1\oplus \cdots \oplus N^f_{\tau_f}
$$
where $N^f_0=f_*\D^\perp$. Setting
$$
\tau^o_f = \begin{cases}
\tau_f\;\;\;\;\;\;\;\;\;\;\;\mbox{if}\;\;n-m\;\;\; \mbox{is even}\\
\tau_f-1\;\;\;\;\; \mbox{if}\;\;n-m\;\;\; \mbox{is odd}
\end{cases}
$$ 
it turns out that the almost complex structure $J$ on $\D^\perp$ induces an almost
complex  structure $J_\ell$ on each $N_\ell^f$, $0\leq \ell\leq\tau^o_f$, defined by
$$
J_\ell\alpha^{\ell+1}_f(X_1,\ldots,X_\ell,X_{\ell+1})
=\alpha^{\ell+1}_f( X_1,\ldots,X_\ell,J X_{\ell+1})
$$
where $\alpha^1_f=f_*$.

The \emph{$\ell^{th}$-order curvature ellipse} 
$\mathcal{E}_\ell^f(x)\subset N^f_\ell(x)$ of $f$ at $x\in M^m$ for 
$0\leq\ell\leq\tau^o_f$ is 
$$
\mathcal{E}_\ell^f(x)=\big\{\alpha_f^{\ell+1}(Z_{\theta},\dots,Z_{\theta}) : 
Z_{\theta}=\cos\theta Z+\sin\theta J Z\;\;\mbox{and}\;\;\theta\in [0,\pi)\big\}
$$
where $Z\in \D^\perp(x)$ has  unit length and satisfies $\<Z,JZ\>=0$. From 
ellipticity such a $Z$ always exists and $\mathcal{E}_\ell^f(x)$ is indeed  
an ellipse. 
\vspace{1ex}

We say that the curvature ellipse $\mathcal{E}_\ell^f$ of an elliptic submanifold
$f$ is a \emph{circle} for some $0\leq\ell\leq\tau^o_f$ if all ellipses 
$\mathcal{E}_\ell^f(x)$ are circles. That the curvature ellipse 
$\mathcal{E}_\ell^f$ is a circle is equivalent to the almost complex  
structure $J_\ell$ being orthogonal.
Notice that $\mathcal{E}_0^f$ is a circle if and only if $f$ is minimal.
An elliptic submanifold $f$ is called \textit{$\ell$-isotropic} if 
all ellipses of curvature up to order $\ell$ are circles.
Then $f$ is called \emph{isotropic} if the ellipses of curvature 
of any order are circles. 
\vspace{1ex}

Substantial isotropic surfaces in $\R^{2n}$ are holomorphic 
curves in $\C^n\equiv\R^{2n}$. Isotropic surfaces in spheres are also 
referred to as  \emph{pseudoholomorphic} surfaces. 
For this class of surfaces a Weierstrass type representation was given 
in \cite{DV}.
\vspace{1ex}

Let  $f\colon M^m\to\Q_c^{n-c}$, $(c=0,1),$ be a substantial nicely curved 
elliptic submanifold. Assume that $M^m$ is the saturation  of a fixed simply 
connected cross section $L^2\subset M^m$ to the relative nullity foliation.
The subbundles in the orthogonal splitting (\ref{splits}) are parallel in 
the normal connection (and thus in $\Q_c^{n-c}$) along $\D$. 
Hence each $N^f_\ell$ can be seen as a vector bundle along the 
surface $L^2$.
\vspace{1ex}

\noindent A \emph{polar surface} to $f$ is an immersion $h$ of $L^2$
defined as follows:
\begin{itemize}
\item [(a)] If $n-c-m$ is odd, then  the polar surface $h\colon L^2\to\Sf^{n-1}$ 
is the spherical image of the unit normal field spanning $N^f_{\tau_f}$. 
\item [(b)] If $n-c-m$ is even, then  the polar surface $h\colon L^2\to\R^n$ is 
any surface such that $h_*T_xL=N^f_{\tau_f}(x)$ up to parallel 
identification in $\R^n$.
\end{itemize}

Polar surfaces always exist since in case $\rm(b)$  any elliptic submanifold 
admits locally many polar surfaces. 
\vspace{1ex}

The almost complex structure $J$ on  $\D^\perp$ induces an almost complex 
structure $\tilde J$ on $TL$ defined by $P\tilde J=JP$
where $P\colon TL \to \D^\perp$ is the orthogonal projection.
It turns out that a polar surface of an elliptic submanifold is necessarily elliptic.  
Moreover, if the elliptic submanifold has a circular ellipse of curvature then its 
polar surface has the same property at the ``corresponding" normal bundle.
As a matter of fact, up to parallel identification it holds that
\be\label{eqp}
N_s^h=N_{\tau^o_f-s}^f\;\;\mbox{and}\;\;
J^h_s=\big(J^f_{\tau^o_f-s}\big)^t,\;\; 0\leq s\leq\tau^o_f.
\ee
In particular, the polar surface is nicely curved.
Notice that the last $\ell+1$ ellipses  of curvature of the  polar surface 
to an $\ell$-isotropic submanifold  are circles. Note that in this case 
the polar surface is not necessarily minimal.
\vspace{1ex}

A \emph{bipolar surface} to $f$  is any polar surface to a polar surface to $f$.
In particular, if we are in case $f\colon M^3\to\Sf^{n-1}$, then a bipolar
surface to $f$ is a nicely curved elliptic immersion $g\colon
L^2\to\R^n$.

\section{The local case}\label{polar-bip}

We discuss next two alternative ways to parametrically describe, at least 
locally, all spherical three-dimensional minimal submanifolds of rank
two in spheres. This follows from the results in \cite{DF2} bearing
in mind that a submanifold is minimal in a sphere if and only if the cone 
shaped over it is minimal in the Euclidean space.
\vspace{1ex}

Let $g\colon L^2\to\R^{n+1}$, $n\ge 4$, be an elliptic surface and let 
$T^1L$ denote its unit tangent bundle.

\begin{proposition}\po\label{local-3dim}
If $\mathcal{E}_1^g$ is a circle, then the map $\psi_g\colon T^1L\to\Sf^n$ 
given by  
$$
\psi_g(x,w)=g_*w
$$
is a minimal immersion with index of relative nullity $\nu\geq 1$ outside the 
subset of singular points, which correspond to points where 
$\dim N_1^g=0$. Moreover, a regular point \mbox{$(x,w)\in T^1L$} is totally 
geodesic for $\psi_g$ if and only if $\dim N_2^g(x)=0$.
Conversely, any three-dimensional minimal submanifold in the sphere with $\nu=1$ 
at any point can be at least locally parametrized in this way.
\end{proposition}

The above parametrization (used for Theorem \ref{main}) is called  
the \emph{bipolar parametrization} in \cite{DF2} because $g$ is a bipolar surface to 
$\psi_g$. The parametrization in the sequel (used for the examples discussed above) 
was called in \cite{DF2} the  \emph{polar parametrization}. 
\vspace{1ex}

Let $g\colon L^2\to\Q_{1-c}^{2n+2c}$  $(c=0,1)$, $n\ge 2$, be a nicely 
curved elliptic surface and let $M^3=UN^g_{\tau_g}$ stand for the unit 
subbundle of $N^g_{\tau_g}$.

\begin{proposition}\po\label{polar-cl}
If $\mathcal{E}^g_{\tau_g-1}$ is a circle, then 
$\phi_g\colon M^3\to\Sf^{2n+c}$ given by 
$\phi_g(x,w)=w$ is a minimal immersion of rank two and polar surface $g$. 
Conversely, any  minimal submanifold $M^3$ in $\Sf^{2n+c}$ of 
rank two can locally be parametrized in this way.
\end{proposition}

\section{The complete case}

We first observe that for complete submanifolds of rank at most two the 
interesting case is the three-dimensional one.  The remaining of the paper
is devoted to the study of the latter case.

\begin{proposition}\label{3d}\po
Let $f\colon M^m\to\Sf^n$, $m\geq 3$, be a minimal isometric immersion with 
index of relative nullity $\nu\geq m-2$ at any point. If $M^m$ is complete, 
then $f$ is totally geodesic unless $m=3$.
\end{proposition}

The above is an immediate consequence of the following result due to Ferus \cite{fe}
(see \cite[Lemma 6.16]{da} where the proof  holds regardless the codimension) since 
due to minimality we cannot have points with index of relative nullity $m-1$.

\begin{lemma}\label{Cr}\po Let $f\colon M^m\to\Sf^n$ be an isometric immersion
and let $U\subset M^m$ be an open subset where the index of relative nullity is 
constant either $\nu=m-1$ or $\nu=m-2$.  Then no leaf of relative nullity contained in 
$U$ is complete if $m\geq 4$.  
\end{lemma}

In the sequel, let $f\colon M^3\to\Sf^n$ be a minimal isometric immersion of a 
complete Riemannian manifold with index of relative nullity $\nu(x)\geq 1$ 
at any $x\in M^3$. Let $U\subset M^3$ be an open subset where  $\nu=1$ such that 
$\D$ is a line bundle on $U$. If that line bundle is trivial, then there is a unique, 
up to a sign,  orthogonal almost complex structure $J\colon\D^\perp|_U\to\D^\perp|_U$. 
In that case set $\mathcal{C}=C_e$ where $e$ is a unit section of $\D|_U$. 

\begin{lemma}\label{h}\po
If  $\D|_U$ is a trivial line bundle  there are harmonic functions 
$u,v\in C^{\infty}(U)$ such that
\be\label{C}
\mathcal{C}=vI-uJ.
\ee
\end{lemma}

\proof
The proof follows similarly as in \cite{dksv}.
Without loss of generality we assume that $f$ is substantial.
Denote by $A_\xi$ the shape operator of $f$ with respect to $\xi\in N_fM$.
The Codazzi equation gives
$$
\nabla_{e}A_\xi|_{\D^{\perp}}= A_\xi|_{\D^{\perp}}
\circ\mathcal{C}+A_{\nabla^{\perp}_e\xi}|_{\D^{\perp}}
$$
for any vector field $\xi\in N_fM$. In particular,
\be\label{eq1}
A_\xi|_{\D^{\perp}}\circ \mathcal{C}=\mathcal{C}^t\circ A_\xi|_{\D^{\perp}}.
\ee
Moreover, the minimality condition is equivalent to
\be\label{eq2}
A_\xi|_{\D^{\perp}}\circ J =J^t\circ A_\xi|_{\D^{\perp}}.
\ee

First we consider the case $n=4$. Let $e_1,e_2,e_3=e$ be a local 
orthonormal
frame field that diagonalizes $A_\xi$   with respect to a unit normal
vector filed $\xi$ such that $J e_1=e_2$. Set
$$
u=\<\n_{e_2}e_1,e_3\>\quad\text{and}\quad v=\<\n_{e_1}e_1,e_3\>.
$$
From the Codazzi equations
$$
(\n_{e_i} A_\xi)e_3=(\n {e_3}A_\xi)e_i\;\; 1\leq i\leq 2,
$$
we obtain that $\<\n_{e_2}e_2,e_3\>=v$ and from the Codazzi equation
$$
\<(\n_{e_1}A_\xi)e_2,e_3\> =\<(\n_{e_2}A_\xi)e_1,e_3\>
$$
that $\<\n_{e_1}e_2,e_3\> =-u$, and now (\ref{C}) follows.

Assume now that $f$ does not reduce codimension to one. Due to the minimality 
assumption, we have that $\dim N_1^f\le 2$. 
If $\dim N_1^f=1$ on an open subset $V\subset M^3$, 
a simple argument using the Codazzi equation gives that $N_1^f$ is parallel in 
the normal bundle, and hence $f|_V$ reduces codimension to one. 
Due to real analyticity the same would hold globally, and that has been excluded. 
Hence, there is an open dense subset $W\subset M^3$ where 
$\dim N_1^f=2$. From (\ref{eq1}) and (\ref{eq2}) we have that
$\mathcal{C}\in \operatorname{span}\{I,J\}$
on $U\cap W$, and (\ref{C}) follows easily.

We now show that $u$ and $v$ are harmonic functions. From 
(\ref{C1}) and (\ref{C2}) we have 
\be\label{c1}
\nabla^{h}_{e} \mathcal{C}=\mathcal{C}^2+I
\ee
and
\be\label{c2}
\big(\nabla^{h}_{X}\mathcal{C}\big)Y=\big(\nabla^{h}_{Y}\mathcal{C}\big)X
\ee
for any  $X,Y\in\D^{\perp}$.
Let $e_1,e_2,e_3$ be a local orthonormal frame with $Je_1=e_2$ and 
$e_3\in\D$. From  
(\ref{C}) we find that
\be\label{o1}
v=\<\nabla_{e_1} e_1,e_3\>=\<\nabla_{e_2}e_2,e_3\>
\ee
and
\be\label{o2}
u=-\<\nabla_{e_1}e_2,e_3\>=\<\nabla_{e_2}e_1,e_3\>.
\ee
Equation  (\ref{c1}) is equivalent to 
\be\label{v}
e_3(v)=v^2-u^2+1\;\;\mbox{and}\;\; e_3(u)=2uv
\ee
whereas (\ref{c2}) is equivalent to
\be\label{u}
e_1(u)=e_2(v)\;\;\mbox{and}\;\; e_2(u)=-e_1(v).
\ee
The Laplacian of  $v$ is given by
$$
\Delta v=\sum_{i=1}^3e_ie_i(v)
-\omega_{12}(e_1)e_2(v)-\omega_{13}(e_1)e_3(v)
+\omega_{12}(e_2)e_1(v)-\omega_{23}(e_2)e_3(v)
$$
where $\omega_{ij}(e_k)=\<\nabla_{e_k}e_i,e_j\>$ for any $i,j,k\in\{1,2,3\}$. 
Using (\ref{u}) we have
\bea
e_1e_1(v)+e_2e_2(v)\!\!\!&=&\!\!\!-e_1e_2(u)+e_2e_1(u)
=-[e_1,e_2](u)
=-\nabla_{e_1}e_2(u)+\nabla_{e_2}e_1(u)\\
\!\!\!&=&\!\!\!\omega_{12}(e_1)e_1(u)-\omega_{23}(e_1)e_3(u)
+\omega_{12}(e_2)e_2(u)+\omega_{13}(e_2)e_3(u)\\
\!\!\!&=&\!\!\!\omega_{12}(e_1)e_2(v)-\omega_{12}(e_2)e_1(v)+2ue_3(u).
\eea
Inserting this equality into the previous equation and making 
use of  (\ref{o1}) and  (\ref{v}) yields $\Delta v=0$.
In a similar form it follows that also $u$ is harmonic.
\vspace{1,5ex}\qed

Let  $\mathcal{A}$ denote the set of totally geodesic points of $f$. 
By Proposition \ref{comp} the relative nullity distribution $\D$ 
is a line bundle on $M^3\smallsetminus\mathcal{A}$.
Being $f$ real analytic, the square of the norm of 
the second fundamental form is a real analytic function and hence 
$\mathcal{A}$ is a real analytic set.  According to Lojasewicz's structure 
theorem \cite[Theorem 6.3.3]{kr} the set $\mathcal{A}$ locally decomposes as
$$
\mathcal{A}=\mathcal{V}^0\cup\mathcal{V}^1\cup \mathcal{V}^2\cup\mathcal{V}^3
$$
where each $\mathcal{V}^d,\, 0\leq d\leq3$, is either empty or a 
disjoint finite union of $d$-dimensional real analytic subvarieties. 
A point $x_0 \in\mathcal{A}$ is called a \emph{regular point of dimension} $d$ 
if there is a neighborhood $\Omega$ of $x_0$ such that $\Omega\cap\mathcal{A}$  
is a $d$-dimensional real analytic submanifold of $\Omega$. 
Otherwise $x_0$ is said to be a \emph{singular} point.
The set of singular points is locally a finite union of submanifolds.

  We want to show that $\mathcal{A}=\mathcal{V}^1$ unless $f$ is
just a totally geodesic three-sphere in $\Sf^n$. 
After excluding the latter case, we have from the real analyticity of $f$
that $\mathcal{V}^3$ is empty. We will proceed now following ideas developed 
in \cite{dksv}. 
In fact, we only sketch the proof of the following fact, which is similar 
to the proof of Lemma $2$ in \cite{dksv}.

\begin{lemma}\label{gp2}\po
The set $\mathcal{V}^2$ is empty.
\end{lemma}

\proof We only have to show that there are no regular points in $\mathcal{V}^2$.
Suppose that a regular point  $x_0\in\mathcal{V}^2$ exists. Let $\Omega \subset M^3$ 
 be an open neighborhood of $x_0$ such that $L^2=\Omega\cap\mathcal{A}$ is an embedded 
surface. Let $e_1,e_2,e_3,\xi_1,...,\xi_{n-3}$ be an orthonormal frame  adapted to $M^3$
along $\Omega$ near $x_0$. The coefficients of the second fundamental form are
$h^a_{ij}=\<\a_f(e_i,e_j),\xi_a\>$ where $1\leq i,j,k\leq 3$ and $1\leq a,b\leq n-3$.

The Gauss map $\gamma\colon M^3 \to Gr(4,n+1)$ of $f$ is a map into the Grassmannian of 
oriented  $4$-dimensional subspaces in $\R^{n+1}$ defined by 
$$
\gamma=f\wedge e_1\wedge e_2 \wedge e_3
$$
where we regard $Gr(4,n+1)$ as a  submanifold in $\wedge^4\R^{n+1}$ via the map for the
Pl\"{u}cker embedding. Then
$$
\gamma_*e_i=\sum_{j,a} h^a_{ij}f\wedge e_{ja}
$$
where $ e_{ja}$ is taken by replacing $e_j$ with $\xi_a$ in $e_1\wedge e_2 \wedge e_2$.
Moreover, it easy to  see that the Gauss map satisfies the partial 
differential equation
$$
\Delta\gamma+\|\a_f\|^2\gamma
=\sum_{i,a\neq b,j\neq k}h^a_{ij}h^b_{ik}f\wedge e_{ja,kb}
$$
where $e_{ja,kb}$ is obtained by replacing $e_j$ with $\xi_a$ and $e_k$ 
with $\xi_b$ in $e_1\wedge e_2\wedge e_3$. Hence, we may write the 
latter equation in the form
$$
\Delta \gamma(x)+\|\gamma_*(x)\|^2\gamma(x)+G(x,\gamma_*)=0
$$
where $G$ is real analytic with $G(\cdot\,,0)=0$. 
Clearly, we have that $\gamma$ is constant along 
$L^2$ and that $\gamma_*({\bf n})=0$ on $L^2$, where ${\bf n}$ is a unit normal of  
$L^2\subset M^3$. Then, it follows from the uniqueness part of the Cauchy-Kowalewsky 
theorem (cf.\ \cite{t}) that the Gauss map $\gamma$ must be constant.
This would imply that $f(M)$ is a three-dimensional totally geodesic sphere 
which contradicts our assumption. \qed

\begin{lemma}\label{gp3}\po
The set $\mathcal{V}^0$ is empty.
\end{lemma}

\proof Let $\Omega$ be an open neighborhood around  $x_0\in\mathcal{V}^0$ such that 
$\nu=1$ on  $\Omega\smallsetminus\{x_0\}$ and let $\{x_j\}_{j\in\mathbb{N}}$ be a 
sequence in $\Omega\smallsetminus\{x_0\}$ converging to $x_0$. 
Let $e_j=e(x_j)\in T_{x_j} M$ be the sequence of unit vectors contained in the 
relative nullity distribution of $f$. 
By passing to a subsequence, if necessary, there is a unit vector 
$e_0\in T_{x_0} M$ such that $\lim e_j=e_0$. By continuity, the geodesic
tangent to $e_0$ at $x_0$ is a leaf of relative nullity outside $x_0$.
But this is a contradiction in view of Proposition \ref{comp}.\qed

\begin{lemma}\label{gp5}\po
The foliation of relative nullity extends analytically over 
the regular points in the set $\mathcal{A}$.
\end{lemma}

\proof Because $\mathcal{A}=\mathcal{V}^1$ its $2$-capacity 
${\mathrm {cap}}_2(\mathcal{A})$ must be zero (cf.\ \cite[Theorem 3]{ev}). 
On the other hand, the distribution $\D$ 
extends continuously over the regular points of $\mathcal{A}$.
In fact, by the previous lemmas  it remains to consider the case when 
$\Omega$ is an open subset of $M^3$ such that $\Omega\cap\mathcal{A}$  
is a open piece of a great circle in the ambient space. 
But in  this situation the result follows by an argument of continuity 
similar than in the proof of Lemma \ref{gp3}. 

Let $\Omega$ be an open subset of $M^3\smallsetminus\mathcal{A}$ 
and $e_1,e_2,e_3$ a local frame on $\Omega$ as in the proof of Lemma \ref{h}. 
Consider the map  $F\colon\Omega\to\Sf^n$ with values into the unit 
sphere given by $F=f_*e_3$. A straightforward computation using 
(\ref{o1}), (\ref{o2}) and (\ref{u})  shows that its tension field
$$
\tau(F)=\sum_{j=1}^3\left(\bar\nabla_{F_*e_j}F_*e_j-F_*\nabla_{e_j}e_j\right)
$$
vanishes, where $\bar\nabla$ denotes the Levi-Civita connection of $\Sf^n$. 
Hence $F$ is a harmonic map. 
Since $F$ is continuous on $M^3$ and because ${\mathrm {cap}}_2(\mathcal{A})=0$, 
it follows from a result of Meier \cite[Theorem $1$]{me} that $F$ is of class
$C^2$ on $M^3$. But then $F$ is real analytic by a result due to Eells and Sampson 
\cite[Proposition  p.\ 117]{ee}.\qed

\begin{lemma}\label{gp6}\po
The set  $\mathcal{A}$ has no singular points.
\end{lemma}

\proof Let $x_0\in\mathcal{A}$ be a singular point. From 
Lemmas \ref{gp2} and \ref{gp3} the set $\mathcal{A}$ contains 
subvarieties of dimension one. It is well known that the singular points of 
such curves are isolated (cf.\  \cite[Theorem 6.3.3]{kr}). 
Moreover, according to Lemma \ref{gp5} the set of regular points of 
$\mathcal A$ contains geodesic curves of the relative nullity foliation. 
Hence $x_0$ is an intersection of such geodesic curves. Let $\Omega\subset M^3$ 
be an open subset containing $x_0$ such that the restriction of 
$f|_\Omega$ is injective. Consider a fixed cross section $L^2$ to $\D$
on $\Omega\smallsetminus\{x_0\}$. 
Note that the immersion $f$ can be locally parametrized by the embedding 
$\phi\colon L^2\times\Sf^1\to\Sf^n$ given by
$$
\phi(x,t)=\exp_{f(x)}\big(tf_*e\big)=\cos t\,f(x)+\sin t\,f_*e(x)
$$
where $e\in\D|_{L^2}$.
Since $x_0$ is an intersection point of geodesics in the relative nullity 
foliation, it follows from the parametrization that there are
points $(x_1,t_1),(x_2,t_2)\in L^2\times\Sf^1$ such that
$$
\phi(x_1,t_1)=f (x_0)=\phi(x_2,t_2),
$$
which leads to a contradiction. \qed

\section{The proofs}

The proof of our main result relies heavily on the following consequence
of the Omori-Yau maximum principle.

\begin{lemma}\label{maxprinc2}\po
Let $M^m$ be a Riemannian manifold for which the Omori-Yau maximum principle 
holds. If $\varphi\in C^{\infty}(M)$  satisfies the partial 
differential inequality $\Delta\varphi\ge 2\varphi^2$, then $\sup\varphi= 0$. 
In particular, if $\varphi\ge 0$ then $\varphi = 0$.
\end{lemma}

\proof See \cite[Theorem 2.8]{amr} or \cite{hsv2}.\qed
\vspace{2ex}

\noindent\emph{Proof of Theorem \ref{main}:}
By Proposition \ref{3d} we only have to consider 
the case $m=3$. We distinguish two cases.
\vspace{1,5ex}

\noindent\emph{Case $\mathcal{A}=\emptyset$}. 
At first suppose that the line bundle $\D$ is trivial with 
$e$ a unit  global section. By Lemma \ref{h} there are harmonic functions 
$u,v\in C^{\infty}(M)$ such that $\mathcal{C}=vI-uJ.$

We claim that $u$ is nowhere zero. To the contrary suppose that 
$u(x_{0})=0$ at $x_{0}\in M^3$.  Let $\gamma\colon \R\to M^3$ the maximal 
integral curve of $e$ emanating from $x_{0}$.  
The second equation in (\ref{v}) gives that $u$ must vanish along $\gamma$. 
Thus the first equation in (\ref{v}) reduces to
$v'(s)= v^2(s)+1$, where $v(s)= v(\gamma(s))$ is an entire function.  
But this is a contradiction since this equation has no entire solutions. 
In the sequel we  assume that $u>0$. 

Using  (\ref{v}) one can easily see that 
\bea
\Delta\big((u-1)^2+v^2\big)\!\!\!&=&\!\!\! 2(\|\nabla u\|^2+\|\nabla v\|^2)
\geq 2((e(u))^2+(e(v))^2)\\
\!\!\!&\geq&\!\!\! 2\big((u-1)^2+v^2\big)^2,
\eea
and Lemma \ref{maxprinc2} implies that  $\mathcal{C}=-J$.
\vspace{1ex}

Let $\mathcal{U}\subset M^3$ be the open dense subset where 
$f$ is nicely curved. Then let $U\subset\mathcal{U}$ be an open connected subset of  
$\mathcal{U}$ that is the saturation of a simply connected cross section 
$L^2\subset U$ to the relative nullity foliation.  
Hereafter we  work on $U$ where $f$ is  nicely curved. 
Hence polar and bipolar surfaces of $f|_U$ are well defined.

Let $h$ be a polar surface to $f|_U$.  
We have seen that the almost complex structure $J$ on  $\D^\perp$ 
induces an almost complex structure $\tilde J$ on $TL$ defined by $P\tilde J=JP$,
where $P\colon TL\to\D^\perp$ is the orthogonal projection. Moreover,
$h$ is elliptic with respect to $\tilde J$ and (\ref{eqp}) holds.
In addition, it follows from Proposition \ref{polar-cl} 
that $\mathcal{E}^h_{\tau_h-1}$ is a circle.

We claim that the last curvature ellipse $\mathcal{E}_{\tau_h}^h$  of
$h$ is also a circle. In that case  the bipolar surface 
$g\colon L^2\to\R^{n+1}$ to $f$ is $1$-isotropic, and we are done. 
Observe that
$$
N^h_{\tau_h}={\mathrm{span}}\{\xi,\eta\}
$$
where $\xi=f_*e|_{L^2}$ and $\eta=f|_{L^2}$.
Using  $\mathcal{C}=-J$,  we obtain that 
\be\label{J}
\xi_{*}=f_*|_{\D^\perp}\circ J\circ P.
\ee
Consider vector fields $X_1,\dots, X_{\tau_h}, Y\in TL$. Since 
$N^h_{\tau_h-1}=N^f_{0}=f_*(\D^\perp),$ we  have 
$$
\a_h^{\tau_h}(X_1,\dots,X_{\tau_h})=f_*Z
$$
for some $Z\in \D^\perp$. 
Keeping in mind the bundle isometries, we obtain that
$$
\a_h^{\tau_h+1}(X_1,\dots, X_{\tau_h},Y)
=\big(\nabla^{h\perp}_{Y}\a_h^{\tau_h}(X_1,\dots,X_{\tau_h})\big)_{N^h_{\tau_h}}
=\big(\nab_{Y}f_*Z\big)_{N^h_{\tau_h}}.
$$
Taking into account (\ref{J}) we see that
\bea
\a_h^{\tau_h+1}(X_1,\dots, X_{\tau_h},Y)
\!\!\!&=&\!\!\!\<\nab_Yf_*Z, \xi\>\xi
+\<\nab_Yf_*Z,\eta\>\eta\\
\!\!\!&=&\!\!\!-\<f_*Z, \xi_*Y\>\xi-\< f_*Z,\eta_*Y\>\eta\\
\!\!\!&=&\!\!\!-\< Z,  JP Y\>\xi-\< Z,  P Y\>\eta.
\eea
Recall that the almost complex structure $J^h_{\tau_h}$ on $N^h_{\tau_h}$ 
is given by
$$
J^h_{\tau_h}\a_h^{\tau_h+1}(X_1,\dots, X_{\tau_h},Y)
=\a_h^{\tau_h+1}(X_1,\dots,X_{\tau_h}, \tilde J Y).
$$
Since 
$$
\a_h^{\tau_h+1}(X_1,\dots,X_{\tau_h},Y)=-\<Z,JPY\>\xi-\< Z, P Y\>\eta
$$
and
$$
\a_h^{\tau_h+1}(X_1,\dots,X_{\tau_h},\tilde J Y)=\<Z,P Y\>\xi-\< Z,JP Y\>\eta,
$$
we have that the vectors $\a_h^{\tau_h+1}(X_1,\dots,X_{\tau_h},Y)$
and $\a_h^{\tau_h+1}(X_1,\dots,X_{\tau_h},\tilde J Y)$
are perpendicular of the same length. Thus $J^h_{\tau_h}$ is orthogonal,
and proves the claim.

Finally, if the line bundle $\D$ is not trivial, it suffices to argue
for a $2$-fold covering  $\Pi\colon \tilde M^3 \to M^3$
such that the nullity distribution $\tilde \D$ of $ \tilde f =f\circ\Pi$ 
is a trivial line bundle  and $\Pi_*\tilde  \D=\D$.  
\vspace{1ex}

\noindent\emph{Case $\mathcal{A}\neq\emptyset$}. We have seen that the  
relative nullity distribution $\D$  can be extended analytically to a line 
bundle on $M^3$, denoted again by $\D$, over the set of totally geodesic 
points $\mathcal{A}$. 
Without loss of generality, we may assume that there is a global unit section 
$e\in\D$, since otherwise we can pass to the $2$-fold covering space
$$
\tilde M^3 = \{(x,w):x\in M^3,w\in\D(x)\;\mbox{and}\;\|w\|=1\}
$$
and argue as in the previous case. From  Lemma \ref{h}, we know that there  exist 
harmonic functions $u,v\in C^{\infty}(M^3\smallsetminus\mathcal{A})$ such that
(\ref{C})  holds on $M^3\smallsetminus\mathcal{A}$. By previous
results the functions $u$ and $v$ can be extended analytically to 
harmonic functions on the entire $M^3$. Moreover, since $u$ is positive 
on $M^3\smallsetminus\mathcal{A}$ and $\mathcal{A}$
consists of geodesic curves, by continuity we get that $u\ge 0$ on $M^3$. 
Then $\|\mathcal{C}+J\|^2$ is globally well defined and, arguing as in 
the previous case, we conclude again that $\mathcal{C}=-J$ on $M^3$. 
The remaining of the proof now goes as before.\qed
\vspace{1,5ex}

\noindent\emph{Proof of Corollary \ref{Lag}:} By a result of Ejiri \cite{ej}
we have that $f$ is minimal. Let $e_1,e_2,e_3$ be a local 
orthonormal tangent frame such that $e_3\in\D$. Since $f$ is Lagrangian,  
we have that $Je_1,Je_2,Je_3$ is an orthonormal frame in the normal bundle 
of $f$. Moreover, it is well known that the $3$-linear tensor $h$ given by
$$
h(e_i,e_j,e_k)=\<\a_f(e_i,e_j),Je_k\>,\;\;i,j,k\in\{1,2,3\},
$$
is fully symmetric. Away from the totally geodesic points, it is easy 
to check that the vector fields $\a_f(e_1,e_1)$ and $\a_f(e_1,e_2)$ are 
perpendicular to each other and have the same length. 
Hence $\mathcal{E}^f_1$ is a circle. 

Suppose at first that $f$ is substantial in $\Sf^6$. Assume that 
the submanifold is the saturation of a fixed cross section $L^2$ 
to the relative nullity foliation. Let $h\colon L^2\to\Sf^6$ be a 
polar surface  to  $f$. From (\ref{eqp}) we obtain that $h$ is 
$1$-isotropic. Proceeding as in the proof of Theorem \ref{main}, 
we deduce that the second ellipse of $h$ is also a circle. 
Therefore, $h$ is pseudoholomorphic and any bipolar surface $g$ to $f$ 
is $2$-isotropic in $\R^7$.

Now we consider the case where $f$ is substantial in $\Sf^5$. 
Take a fixed cross section $L^2$ to the relative nullity foliation 
and let $h\colon L^2\to\R^6$ be a polar surface  to $f$. 
As in the previous case, we obtain that $h$ must be isotropic. 
Therefore, any bipolar surface $g$ to $f$ is an isotropic surface 
in $\R^6$.\qed

\noindent Marcos Dajczer\\
IMPA -- Estrada Dona Castorina, 110\\
22460--320, Rio de Janeiro -- Brazil\\
e-mail: marcos@impa.br

\bigskip

\noindent Theodoros Kasioumis\\
University of Ioannina \\
Department of Mathematics\\
Ioannina--Greece\\
e-mail: theokasio@gmail.com

\bigskip

\noindent Andreas Savas-Halilaj\\
Leibniz Universit\"at Hannover \\
Institut f\"ur Differentialgeometrie\\
Welfengarten 1\\
30167 Hannover--Germany\\
e-mail: savasha@math.uni-hannover.de

\bigskip

\noindent Theodoros Vlachos\\
University of Ioannina \\
Department of Mathematics\\
Ioannina--Greece\\
e-mail: tvlachos@uoi.gr

\end{document}